\documentclass[10pt,a4paper]{article}
\usepackage[ansinew]{inputenc}
\usepackage{amsmath}
\usepackage{amsfonts}
\usepackage{amssymb}
\usepackage [dvips] {graphicx}
\newtheorem{df}{Definition} 
 
\newtheorem{prop} {Proposition} 
 
\newtheorem{thm} [prop] {Theorem}

\begin{document}

\title{The Connectivity Order of Links}

\author{S. Dugowson\\Institut Sup\'{e}rieur de M\'{e}canique de Paris\\ dugowson@ext.jussieu.fr}

\maketitle

\textsc{Abstract} --- We associate at each link a connectivity space which describes its splittability properties. Then, the notion of order for finite connectivity spaces results in the definition of a new numerical invariant for links, their connectivity order. A section of this short paper presents a theorem which asserts that every finite connectivity structure can be realized by a link : the Brunn-Debrunner-Kanenobu Theorem.

\mbox{}


\mbox{}


\mbox{}

\section{Order of Finite Connectivity Spaces}
 
Let us recall the definition of a connectivity space \cite{Borger:1983,Dugowson:2007c}.

\begin{df} [Connectivity spaces]  A \emph{connectivity space} is a couple $\mathbf{X}=(X,\mathcal{K})$ where $X$ is a set and $\mathcal{K}$ is a set of nonempty subsets of  $X$ such that
\begin{displaymath}
\forall \mathcal{I}\in \mathcal{P}(\mathcal{K}), \bigcap_{K\in\mathcal{I}}K\ne\emptyset\Longrightarrow \bigcup_{K\in\mathcal{I}}K\in\mathcal{K}.
\end{displaymath}
\noindent The set $X$ is called the \emph{support} of the space $\mathbf{X}$, the set $\mathcal{K}$ is its \emph{connectivity structure}. The elements of  $\mathcal{K}$ are called the \emph{connected subsets} of $\mathbf{X}$.
The morphisms between two connectivity spaces are the functions which transform connected subsets into connected subsets.
A connectivity space is called \emph{integral} if every singleton subset is connected. A connectivity space is called \emph{finite} if the number of its points is finite.
\end{df}

\begin{df} [Irreducible connected subsets] Let $\mathbf{X}=(X,\mathcal{K})$ be a connectivity space. A connected subset $K\in\mathcal{K}$ is called \emph{reducible} if there are two connected subsets $A\subsetneqq K$ and $B\subsetneqq K$ such that
\begin{displaymath}
K=A\cup B \textrm{ and }
A\cap B \neq \emptyset.
\end{displaymath}
A connected subset is said to be \emph{irreducible} if it is not reducible.
\end{df}

In the sequel, all connectivity spaces will be supposed to be integral and finite.

\begin{df} [Order of irreducible connected subsets] Let $\mathbf{X}=(X,\mathcal{K})$ be a (finite and integral) connectivity space. We define by induction the \emph{order} of each irreducible subset. Singleton subsets are said to be of order 0. The order $\omega(L)$ of an irreducible subset $L$ which has more than one point is defined by 
\[ 
\omega(L)=1+\max_{K\in \mathcal{S}(L)}\omega(K)
 \]
where $\mathcal{S}(L)$ is the set of irreducible connected subsets which are strictly included in $L$.
\end{df}

\begin{df} [Order of a connectivity space] The order $\omega(\mathbf{X})$ of a (finite and integral) nonempty connectivity space $\mathbf{X}=(X,\mathcal{K})$ is the maximum of the order or its irreducible connected subsets.
\end{df}

Instead of $\omega(\mathbf{X})$, we generally write just $\omega(X)$. 

\mbox{}

\noindent \textbf{Remark}. It is usefull to represent the structure of a (finite and integral) connectivity space $\mathbf{X}$ by a graph $G$ whose vertices are the irreducible connected subsets of $\mathbf{X}$ and edges express the inclusion relations : if $a$ and $b$ are irreducible connected subsets of $\mathbf{X}$, $(a,b)$ is an edge of $G$ iff $a\subseteq b$ and there is no irreducible connected subset $c$ such that $a\subsetneqq c \subsetneqq b$. In \cite{Dugowson:20070717}, I called this graph the \textit{generic graph} of the connectivity space $\mathbf{X}$ ; the vertices of this graph, \textit{i.e.} the irreducible connected subsets of $\mathbf{X}$, are the \textit{generic points} of $\mathbf{X}$. The order of a space is then the maximal length of paths in the generic graph.

\begin{prop} The order $\omega(X)$ of a nonempty finite integral  connectivity space is always less than the number $\sharp X$ of its points : $\omega(X)\leq \sharp X-1$.
The order $\omega(X)$ is zero iff the space is totally disconnected, that is the only connected subsets are the singletons.
\end{prop}

\noindent \textbf{Example}. Let $n\geq 1$ an integer, and $(A_n,\mathcal{A}_n$ the connectivity space defined by $A_n=\{1,\cdots,n\}$ and 
\[ 
\mathcal{A}_n=\{\{1\},\{2\},\cdots,\{n\}\}\cup \{\{1,2\},\{1,2,3\},\cdots,\{1,2,\cdots,n\}\}. 
 \]
Then $\omega(A_n)=n-1$. Indeed, in this space each connected subset is irreducible, and its order is, by induction, equal to its cardinal minus one. It is, up to isomorphism, the only connectivity space with $n$ points which is of order $n-1$.

\section{The Brunn-Debrunner-Kanenobu Theorem}

At each (tame) link $L$, we can associate a connectivity space $\mathbf{X}_L$ taking the components of the link $L$ as points of $\mathbf{X}_L$ and nonsplittable sublinks of $L$ (one considers knots, \textit{i.e.} sublinks with only one component, as nonsplittable links) defining the connected nonempty subsets of $\mathbf{X}_L$.

\begin{df} The connectivity structure of the connectivity space $\mathbf{X}_L$ associated to a link $L$ is called the \emph{splittability structure} of $L$.
\end{df}

\mbox{}

\noindent \textbf{Example}. The splittability structure of the Borromean ring is (isomorphic to) $\{\{1\},\{2\},\{3\},\{1,2,3\}\}$. 

\mbox{}

In \cite{Dugowson:2007c, Dugowson:20070717}, I asked whether every finite connectivity space can be represented by a link, that is whether exists (at least) one link whose connectivity structure is (isomorphic to) the one given. It appears that in 1892, Brunn \cite{Brunn:1892a} first asked this question, without formalizing the concept of a connectivity space. His answer was positive, and he gave the idea of a proof based on a construction using some of the links now called ``brunnian''. In 1964, Debrunner \cite{Debrunner:1964}, rejecting the Brunnian ``proof", gave another construction, proving it but only for $n$-dimensional links with $n\geq 2$. In 1985, Kanenobu \cite{Kanenobu:198504,Kanenobu:1986} seems to be the first to give a proof of the possibility to represent every finite connectivity structure by a classical link, a result which has not been, up to now, very wellknown. The key idea of those different constructions is already in the Brunn's article : for him, what we call ``Brunnian links" are not so interesting in themselves, but for the constructions they allow to make, that is the representation of \textit{all} finite connectivity strutures by links.

\begin{thm} [Brunn-Debrunner-Kanenobu] Every finite connectivity structure is the splittability structure of a link in $\mathbf{R}^3$.
\end{thm}

\section{The Connectivity Order of Links}

At each link $L$, we associate its connectivity order $\omega(L)$,\textit{ i.e.} the order of the connectivity space associated to $L$.

\mbox{}

\begin{figure} 
\begin{center}
\includegraphics [scale=0.2]{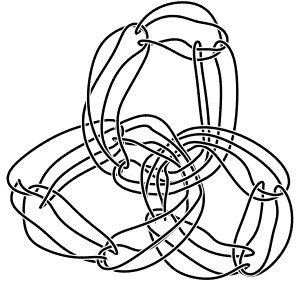}
\caption{A Borromean ring of borromean rings.}
\label{borroborro}
\end{center}
\end{figure}

\begin{figure} 
\begin{center}
\includegraphics [scale=0.3]{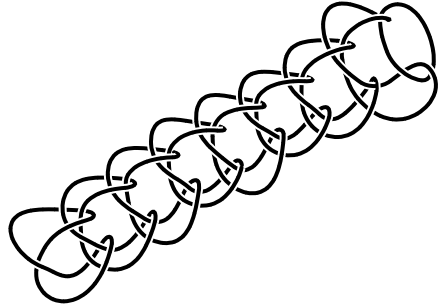}
\caption{A link with a connectivity order 8.}
\label{tendu}
\end{center}
\end{figure}

\noindent \textbf{Examples}. The link pictured on the figure\,\ref{borroborro} has a connectivity order 2. The one of the figure\,\ref{tendu} has a connectivity order 8, which is the maximal order for a link with 9 components. The way this very asymetrical link splits when a component is erased or cut highly depends on the position of this component in the link, as shows its connectivity structure, which is the one we called $\mathcal{A}_9$.

\mbox{}

\noindent \textbf{Remark}. The connectivity order is not a Vassiliev finite type invariant for links. For example, it is easy to check that the connectivity order of the singular link with two components, a circle and another component crossing this circle at $2n$ double-points, is greater than $2^n$.

\bibliographystyle{amsplain}

\end{document}